\documentclass[12pt]{article}

\usepackage{amsmath, amsthm,amsfonts,amssymb,amscd}
\usepackage[dvips]{graphicx}
\usepackage{epic,eepic}
\theoremstyle{plain}
\newtheorem{theorem}{Theorem }[section]
\newtheorem{proposition}[theorem]{Proposition}
\newtheorem{lemma}[theorem]{Lemma}
\newtheorem{corollary}[theorem]{Corollary}
\newtheorem{maintheorem}{Theorem}

\theoremstyle{definition}
\newtheorem{remark}[theorem]{Remark}

\newtheorem{definition}[theorem]{Definition}

\newtheorem{claim}{Claim}

\newtheorem*{ack}{Acknowledgements}
\theoremstyle{remark}

\newcommand{\field}[1]{\mathbb{#1}}
\newcommand{\real}{\field{R}}

\newcommand{\al} {\alpha}       
        
\newcommand{\ga} {\gamma}    
\newcommand{\de} {\delta}       
\newcommand{\ep} {\epsilon}

\newcommand{\la} {\lambda}

\newcommand{\vfi}{\varphi}
\newcommand{\om} {\omega}

\newcommand{\SG}{{\cal G}}

\newcommand{\SL}{{\cal L}}

\newcommand{\SO}{{\cal O}}

\renewcommand{\SS}{{\cal S}}

\newcommand{\SU}{{\cal U}}
\newcommand{\SV}{{\cal V}}

\newcommand{\eps}{\varepsilon}

\newcommand{\mundo}{\operatorname{Diff}^1_{\omega}(M)}
\newcommand{\submundo}{\operatorname{Diff}^{1+\al}_{\omega}(M)}
\newcommand{\newmundo}{\operatorname{Diff}^{k+\al}_{\omega}(M)}

\newcommand{\Mundo}{\mathfrak{X}^{1}_m(M)}
\newcommand{\Submundo}{\mathfrak{X}^{\infty}_m(M)}
\newcommand{\Subsubmundo}{\mathfrak{X}^{1+\al}_m(M)}
\newcommand{\Newmundo}{\mathfrak{X}^{k+\al}_{\omega}(M)}

\newcommand{\N}{\mathbb{N}}
\newcommand{\R}{\mathbb{R}}

\newcommand{\ov}{\overline}

\begin{document}

\title{A pasting lemma and some applications for conservative systems}
\author{Alexander Arbieto and Carlos Matheus}
\date{January 10, 2006}

\maketitle

\begin{abstract}
We prove that in a compact manifold of dimension $n\geq 2$, a
$C^{1+\alpha}$ volume-preserving diffeomorphisms that are robustly
transitive in the $C^1$-topology have a dominated splitting. Also we
prove that for 3-dimensional compact manifolds, an isolated robustly
transitive invariant set for a divergence-free vector field can not
have a singularity. In particular, we prove that robustly transitive
divergence-free vector fields in 3-dimensional manifolds are Anosov.
For this, we prove some ``pasting'' lemma, which allows to make
perturbations in conservative systems.
\end{abstract}

%%%%%%%%%%%%%%%%%%%%%%%%%%%%%%%%%%%%%%%%%%%%%%%%%

\section{Introduction}
\label{intro}

%This is a preliminary version of a work in progress. So, we have to explain why
%we are writing it. The explanation is somewhat simple: in the present work we
%proved some results concerning some ``pasting'' (perturbation) lemma which
%were used by another researchers, although there is not even a \emph{preprint}
%version of the theorems. However, the disadvantage of this preliminary version
%is that it is too close to a draft of the theorems that we hope to prove in a
%(if possible) near by future. But, we point out that the results below are
%correct, although some of them need a refinement. In particular, the suggestions
%of the readers of this ``pre-preprint'' certainly will help in order to improve
%this material.

In this article, we study some properties of conservative dynamical
systems over a smooth $n$-dimensional Riemannian manifold $M$
without boundary, in the discrete and continuous case (here $n\geq
2$). The main technical result in this work is a general tool for
perturbing $C^2$ conservative dynamical systems over (large) compact
subsets of phase space, obtaining a globally defined and still
conservative perturbations of the original system. In fact, we give
several versions of these pasting lemmas, for discrete as well as
for continuous time systems (and in this case is enough suppose that
the system is $C^1$). Before hand, let us describe some of its
consequences in dynamics. In the rest of the paper, the Riemannian
metric is $C^{\infty}$.

\smallskip
{\it Robustly transitive diffeomorphisms}

A dynamical system is \emph{robustly transitive} if any $C^1$ nearby one has
orbits that are dense in the whole ambient space. Recall that a \emph{dominated splitting}
is a continuous decomposition
$TM=E\oplus F$ of the tangent bundle into continuous subbundles which are
invariant under the derivative $Df$ and such that $Df|_E$ is more expanding/less
contracting than $Df|_F$ by a definite factor (see precise definitions and more
history on section~\ref{s.dom}).

We prove: ``\emph{Let $f:M^n\to M^n$ a $C^{1+\alpha}$ volume
preserving diffeomorphism robustly transitive among volume
preserving diffeomorphisms, where $n\geq 2$. Then $f$ admits a
dominated splitting of the tangent bundle}''. This is an extension
of a theorem by Bonatti, Diaz and Pujals.

Also this solves a question posed by Tahzibi~\cite{T}: every $C^1$-stably
ergodic $C^{1+\eps}$-diffeomorphism has a dominated splitting. We observe that
this result was used by Bochi, Fayad and Pujals~\cite{BFP} to prove that there exists an open and
dense set of $C^1$-stably ergodic $C^{1+\eps}$-diffeomorphisms that are
non-uniformly hyperbolic.

\smallskip

{\it Robustly transitive vector fields}

First we prove an extension of a theorem by Doering in the conservative setting:
``\emph{Let $X$ be a divergent free vector field robustly transitive among
divergence free vector fields in a three dimensional manifold. Then $X$ is Anosov.}

It has been shown by Morales-Pac\'{\i}fico and Pujals~\cite{MPP} that robustly
transitive sets $\Lambda$ of a 3-dimensional flows are \emph{singular
hyperbolic}: there exists a dominated splitting of the tangent bundle restricted to $\Lambda$ such
that the one-dimensional subbundle is hyperbolic(contracting or expanding and
the complementary one is volume hyperbolic). Robust transitivity means that
$\Lambda$ is the maximal invariant set in a neighborhood $U$ and, for any $C^1$
nearby vector field, the maximal invariant set inside $U$ contains dense orbits.

For generic dissipative flows robust transitivity can not imply
hyperbolicity, in view of the phenomena of Lorenz-like attractors
containing both equilibria and regular orbits. However we also prove
that Lorenz-like sets do not exist for conservative flows:
``\emph{If $\Lambda$ is a robust transitive set within divergence
free vector fields then $\Lambda$ contains no equilibrium points.}

\smallskip

The paper is organized as follows. In section~\ref{perturbations} we
prove a result of denseness of volume preserve systems with higher
differentiability among $C^1$ volume preserving systems. In
section~\ref{lemmas} we prove the pasting lemmas for volume
preserving diffeomorphisms and vector fields, and for symplectic
diffeomorphisms. In section~\ref{mpp}, we use the pasting lemma to
study robustly transitive volume preserving flows in 3-manifolds. In
section~\ref{s.dom}, we use the pasting lemma to study robustly
transitive conservative diffeomorphisms. Finally in the appendix, we
investigate why some techniques involved in the proof of the pasting
lemma for diffeomorphisms don't work for $C^1$ diffeomorphisms.
Indeed, we given a ``non-estimate'' argument to get this.

\subsection{Outline of the proofs}

To prove the statement for diffeomorphisms, we fix $f$ a robuslty
transitive volume-preserving diffeomorphism then we use a result of
Bonatti-Diaz-Pujals, that says that either any homoclinic class
$H(p,f)$ has a dominated splitting or there exists a sequence of
volume preserving diffeomorphisms $g_n\to f$ and periodic points
$x_n$ of $g_n$, such that $Dg_n^{p(x_n)}(x_n)=Id$ where $p(x_n)$ is
the period of $x_n$. Now we use the pasting lemma, to obtain another
\emph{volume preserving} $h$ close to $f$ such that in a
neighborhood $U$ of $x_n$ we have $h=id$. This contradicts the
robust transitivity of $f$. An important observation here is that to
use the pasting lemma, we the original diffeomorphism need to be
$C^{1+\eps}$.

Now to prove the statement for vector fields, we follow the arguments of
Morales-Pacifico-Pujals. The idea is that if there exists a singularity, then
this singularity is Lorenz-like, to obtain that we use the pasting lemma to get
a contradiction with the robust transitivity if the singularity is not
Lorenz-like. Then an abstract lemma says that this set is a
proper attractor. So, either the set is the whole manifold or we have a proper
attractor. In the first case we prove that the vector field is Anosov, so there
are not singularities, in the second we have a contradiction because the vector
field is divergence-free. Another point is
that for vector fields, we can use the pasting lemma for $C^1$ vector fields,
since we know that $C^{\infty}$ divergence-free vector fields are dense in
$C^1$ divergence-free vector fields. This regularity is also used to perform a
Shilnikov bifurcation in a technical lemma.

Finally, the pasting lemma says that if we have a conservative system
(diffeomorphism or flow) in $M$ a smooth manifold, and we have another conservative
system defined in some open set $U$ of $M$, then we can ``cut'' the first system
in some small open set of $U$ and ``paste'' the second system there, and the
result system is again conservative. The idea is to ``glue'' the two systems
using unity partitions, so we have another system that is conservative inside
the domains where the unity partitions are equal to 1. But, of course, we have
an ``annulus'' (a domain with boundary in the manifold) where the system is not
conservative (divergence not equal to zero in the flow case or determinant of the
derivative not equal to one in the diffeomorphism case). Now we will fix this
problem, finding another system in this domain which ``fix'' the
non-conservative part, for instance, in the flow case, if $h(x)$ is the
divergence of the system in the domain, then we search for a vector field with
divergence equal to $-h(x)$, then we sum the two vector fields and we have a
conservative vector field in the domain. But of course, we want that at the
boundary the vector field is zero and that we can extend it as zero to the rest of the
manifold. The problem to find that vector field is then a PDE problem, and we
use PDE arguments to solve it and also obtain smallness for the norm of the
vector field that solves the equation. The procedure is the same for
diffeomorphisms.

\section{Denseness Result}
\label{perturbations}

\begin{definition}
Let $M=M^n$ a Riemannian manifold without boundary with dimension
$n\geq 2$. We say that a vector field $X$ is conservative if div
$X=0$, we denote by $X\in \mathfrak{X}_m(M)$.
\end{definition}

Of course if $X$ is conservative, the flow is volume preserving by the Liouville formula.

Now we proving that $C^{\infty}$ conservative vector fields are $C^1$-dense in
$C^1$ conservative vector fields:

\begin{theorem}
\label{l.dense}
$\Submundo$ is $C^1$-dense on $\Mundo$.
\end{theorem}

\begin{proof}
Let us fix $X\in\Mundo$. Locally we proceed as follows: take a
conservative local charts $(U,\Phi)$ (Moser's theorem~\cite{M}) and
get $\eta_{\eps}$ a $C^{\infty}$ Friedrich's mollifier (i.e. let
$\eta$ a $C^{\infty}$ function with support in the ball $B(0,1)$
such that $\int \eta=1$ and take
$\eta_{\eps}(x)=\eps^{-n}\eta(x/\eps)$). If $X=(X_1,\dots, X_n)$ is
the expression of $X$ in local coordinates, we define
$X_{\eps}=(X_1*\eta_{\eps},\dots X_n*\eta_{\eps})$ a $C^{\infty}$
vector field (where the star operation is the convolution). Now we
note that
$\frac{d}{dx_i}(X_i*\eta_{\eps})=(\frac{d}{dx_i}X_i)*\eta_{\eps}$,
so $\text{div }X_{\eps}=0$. Furthermore $X_{\eps}$ converges to $X$
in the $C^1$-topology.

\smallskip

Now we take local charts $(U_i,\phi_i)$, $i=1,..., m$ as above and
$\xi_i$ a partition of unity subordinated to $U_i$, and open sets
$W_i$ such that $W_i\subset V_i:=supp(\xi_i)\subset U_i$ satisfying
$\phi_i(W_i)=B(0,1/3)$ and $\phi_i(V_i)=B(0,2/3)$, $\xi_i|_{W_i}=1$
and $\Omega:=M\backslash int(W_i)$ is a manifold with boundary
$C^{\infty}$, here $B(0,r)$ is the ball with center 0 and radii $r$.
And we fix $C$ as the constant given by theorem~\ref{t.2dm} (below).

In each of those charts
we get $X_i$ a $C^{\infty}$ conservative vector fields $C^1$-close to $X$ by
theorem~\ref{t.2dm} such that if
we take $Y=\sum_i \xi_iX_i$ a $C^{\infty}$ vector field and we denote $g=\text{div
 }Y$ a $C^{\infty}$ function then:

\begin{itemize}
\item $g$ is $\frac{\eps}{2C}$-$C^1$-close to
0 (where $C$ is the constant in the below theorem). Indeed, we have:

$$g=\sum_i \nabla \xi_i \cdot X_i + \xi_i\cdot \text{div }X_i=\sum_i \nabla \xi_i
\cdot X_i.$$

So, it is suffice to take the $X_i$'s $C^1$-close enough to $X$.
\item $Y$ is $\frac{\eps}{2}$-$C^1$-close to $X$.

\end{itemize}

By the divergence theorem we get that $\int_{\Omega}g =0$, since:

$$\int_{\Omega} \text{div }Y=\int_{\partial \Omega} Y\cdot N=-\sum_i
\int_{\partial W_i} X_i\cdot N =-\sum_i \int_{W_i} \text{div } X_i=0.$$

Now, we state the following results from PDE's by
Dacorogna-Moser~\cite[Theorem 2]{DM} on the divergence equation which will
play also a crucial role in the next section in the proofs of the pasting lemmas:

\begin{theorem}
\label{t.2dm}
Let $\Omega$ a manifold with $C^{\infty}$ boundary. Let $g\in
C^{k+\alpha}(\Omega)$ (with $k+\alpha>0$) such that $\int_{\Omega}g=0$. Then there exist $v$ a
$C^{k+1+\alpha}$ vector field (with the same regularity at the boundary) such
that:
\begin{displaymath}
\left\{ \begin{array}{ll}
\text{div }v(x)=g(x) \text{, }x\in\Omega \\
v(x)=0 \text{, }x\in \partial \Omega \end{array}\right.
\end{displaymath}
Furthermore there exist $C=C(\alpha,k,\Omega)>0$ such that $\|v\|_{k+1+\alpha}\leq
C\|g\|_{k+\alpha}$. Also if $g$ is $C^{\infty}$ then $v$ is $C^{\infty}$.
\end{theorem}

\begin{remark}
In fact Dacorogna-Moser prove this result for open sets in $\R^n$. But the same
methods work for manifolds, because it follows from the
solvability of $\triangle u=f$ with Neumann's condition. Also, if the manifold
has boundary, we obtain
solutions that have the required regularity at the boundary (see~\cite[p.3]{DM}
or~\cite[p.265]{H}; see
also~\cite[Ch 4,Theorem 4.8]{Au},~\cite{GT} and~\cite{LU}).
\end{remark}

So, we get $v$ the vector field given by theorem~\ref{t.2dm} and we define
$Z=Y-v$. We observe that $Z$ is a $C^{\infty}$ vector field because $v$ is
$C^{\infty}$ at the boundary. Also $\text{div }Z=0$. Finally since $Y$ is
$\frac{\eps}{2}$-$C^1$-close to $X$ and $v$ is $\frac{\eps}{2}$-$C^1$-close to 0
we get $Z$ $\eps$-$C^1$-close to $X$.
The proof is complete.
\end{proof}

\begin{remark}As Ali Tahzibi pointed out to the second author after the
conclusion of this work, the previous
theorem was proved by Zuppa~\cite{Zu} (in 1979) with a different proof. In fact,
since Dacorogna-Moser's result was not available at that time, Zuppa uses that
the Laplacian operator admits a right inverse.
\end{remark}

\section{The Pasting Lemmas}
\label{lemmas}

\subsection{Vector Fields}

We start proving a weak version of the pasting lemma:

\begin{theorem}[The $C^{1+\alpha}$-Pasting Lemma for Vector Fields]
\label{t.c2glue} Let $M^n$ a compact Riemannian manifold without
boundary with dimension $n\geq 2$. Given $\al>0$ and $\eps_0>0$
there exists $\delta_0>0$ such that if $X_0\in \Subsubmundo$, $K$ is
a compact subset of $M$ and $Y_0\in \Subsubmundo$ is
$\delta_0$-$C^1$-close to $X_0$ on a small neighborhood $U$ of $K$.
Then there exist a $Z_0\in \Subsubmundo$, and $V$ and $W$ such that
$K\subset V\subset U\subset W$ satisfying $Z_0|_V=Y_0$,
$Z_0|_{int(W^c)}=X_0$ and $Z_0$ is $\eps_0$-$C^1$-close to $X_0$.
Furthermore, if $X_0$ and $Y_0$ are $C^{\infty}$ then $Z_0$ is also
$C^{\infty}$.
\end{theorem}

\begin{proof}
Let $V$ be a neighborhood of $K$ with $C^{\infty}$ boundary, compactly contained in
$U$ such that $U$ and $int(V^c)$ is a covering of $M$. Let $\xi_1$ and $\xi_2$
be a
partition of unity subordinated to this covering such that $\xi_1|_V=1$ and
there exist $W\subset U^c$ such that $\xi_2|_W=1$. Let $\Omega=M\backslash
V\cup W$ be a (non empty) manifold with $C^{\infty}$ boundary and we fix $C$ as the constant
given by theorem~\ref{t.2dm}.

Now we choose $\delta_0$ such that if $Y_0$ is $C^{1+\al}$ $\delta_0$-$C^1$-close
to $X_0$ then
$T=\xi_1Y_0+\xi_2X_0$ is $\frac{\eps_0}{2}$-$C^1$-close to $X_0$ and $g=\text{div }T$ is $\frac{\eps_0}{2C}$-$C^1$-close to 0. Indeed, we note
that $g=\nabla \xi_1\cdot Y_0+\nabla \xi_2\cdot X_0$ is $C^1$-close to 0 if we get
$Y_0$ sufficiently $C^1$-close to $X_0$. Also note that by the divergence theorem
(see the proof of the theorem~\ref{l.dense} above) $\int_{\Omega}g=0$. Clearly, $g$ is a
$C^{1+\al}$ function.

So we take $v$ the vector field given by the theorem~\ref{t.2dm} and extend it as 0 in the rest of
the manifold. This extension is a
$C^{2+\al}$-vector field (because we have regularity of $v$ at the boundary)
and $Z_0=T-v$ is a $C^{1+\al}$-vector field $\eps_0$-$C^1$-close to
$X_0$. Now $Z_0=Y_0$ in $V$ so $div Z_0=$ in $V$, also in $W$ we have $Z_0=X_0$
so again $div Z_0=0$, finally in $\Omega$ we have $div Z_0=div(T-v)=div T-div v=$ since
$v$ is a solution of the PDE, so $Z_0$
satisfies the statement of the theorem.
\end{proof}

We extend this theorem to the $C^1$ topology using the denseness theorem:

\begin{theorem}[$C^1$-Pasting Lemma for vector fields]
\label{t.c1glue} Let $M^n$ a compact Riemannian manifold without
boundary with dimension $n\geq 2$. Given $\eps>0$ there exists
$\delta>0$ such that if $X\in \Mundo$, $K$ a compact subset of $M$
and $Y\in \Subsubmundo$ is $\delta$-$C^1$-close to $X$ on a small
neighborhood $U$ of $K$, then there exists a $Z\in \Subsubmundo$ and
$V$ such that $K\subset V\subset U$ satisfying $Z|_V=Y$ and $Z$ is
$\eps$-$C^1$-close to $X$. If $Y\in\Submundo$ then $Z$ is also in
$\Submundo$.
\end{theorem}

\begin{proof}
Let $\eps$ and set $\eps_0=\frac{\eps}{2}$ and let $\delta_0$ be given by
theorem~\ref{t.c2glue}. Let $\sigma=\min \{\frac{\eps}{2},
\frac{\delta_0}{2}\}$, and $\delta=\frac{\delta_0}{2}$.

Now get $X_0$ $\sigma$-$C^1$-close to $X$ by theorem~\ref{l.dense}.
So if $Y$ is $\delta$-$C^1$-close to $X$ then $Y$ is
$\delta_0$-$C^1$-close to $X_0$. By theorem~\ref{t.c2glue} we have
$Z$ $\eps_0$-$C^1$-close to $X_0$ satisfying the conclusions of the
theorem and $Z$ is $\eps$-$C^1$-close to $X$.
\end{proof}

Using the same arguments of Dacorogna-Moser, we can produce perturbations in
higher topologies, provided that the original systems are smooth enough.
We will use the following notation, if $r$ is a real number, we denote by $[r]$ its integer part. In
the lemma below $\widetilde{r}$ denotes a real number such that $\widetilde{r}=r$ if $r$
is an integer and $[r]\leq\widetilde{r}<r$ if $r$ is not an integer. So we
obtain:

\begin{theorem}[Pasting Lemma with Higher Differentiability]
\label{l.gluinglemma} Let $M^n$ a compact Riemannian manifold
without boundary with dimension $n\geq 2$. Let $X\in\Newmundo$ ($k$
integer, $0\leq\al <1$ and $k+\al>1$). Fix $K\subset U$ a
neighborhood of a compact set $K$. Let $Y$ be a vector field defined
in $U_1\subset U$, $U_1$ an open set containing $K$. If
$Y\in\Newmundo$ is sufficiently $C^{r}$-close to $f$ for $1\leq
r\leq k+\al$ ($r$ is real) then there exists some $C^{r}$-small
$C^{k+\al}$-perturbation $Z$ of $X$ and an open set $V\subset U$
containing $K$ such that $Z=X$ outside $U$ and $Z=Y$ on $V$.

Furthermore, there exists a constant $C=C(f,\dim(M),K,U,r)$, $\de_0$ and
$\ep_0$ such
that if $Y$ is $(C\cdot\de^{\widetilde{r}}\cdot\ep)$-$C^{r}$ close to $X$,
$\de<\de_0$, $\ep<\ep_0$ then
$Z$ is $\ep$-$C^{r}$ close to $X$ and $\textrm{support}(Z-X)\subset
B_{\de}(K)$ (= the $\de$-neighborhood of $K$).
\end{theorem}

\subsection{A Weak Pasting Lemma for Conservative Diffeomorphisms}

In this subsection we prove a weak version of the pasting lemma for $C^2$
conservative diffeomorphisms which we will apply to the study of dominated
splittings for conservative systems.

First, we note that a similar result as theorem~\ref{t.2dm} holds for diffeomorphisms, as showed by
Dacorogna-Moser~\cite[Theorem 1, Lemma 4]{DM} (see also~\cite{RY}):

\begin{theorem}\label{l.DM4}Let $\Omega$ be a compact manifold with $C^{\infty}$ boundary.
Let $f,g\in
C^{k+\alpha}(\Omega)$ ($k+\al>0$) be such that $f,g>0$. Then there exists a
$C^{k+1+\alpha}$ diffeomorphism $\vfi$ (with the same regularity at the boundary) such
that:
\begin{displaymath}
\left\{ \begin{array}{ll}
g(\vfi(x))\det(D\vfi(x))=\lambda f(x) \text{, }x\in\Omega \\
\vfi(x)=x \text{, }x\in \partial \Omega \end{array}\right.
\end{displaymath}
where $\lambda=\int g/\int f$. Furthermore there exists
$C=C(\alpha,k,\Omega)>0$ such that $\|\vfi-id\|_{k+1+\alpha}\leq
C\|f-g\|_{k+\alpha}$. \footnote{This constant is invariant by
conformal scaling $x\to rx$} Also if $f,g$ are $C^{\infty}$ then
$\vfi$ is $C^{\infty}$.
\end{theorem}

%Let's find $\phi:B(0,2r)-B(0,r)\to B(0,2r)-B(0,r)$ which solves the
%equation for $f$.
%
%Then let $g:B(0,2)-B(0,1)\to R$ such that $g(y)=f(ry)$ and choose
%$\psi(y)$, which solves for $g$. Let $\phi(x)=r\psi(x/r)$
%
%Then $\det D\phi(x)=\det D\psi(x/r)=g(x/r)=f(x)$. So, $\phi$ solves
%the equation for $g$.
%
%Now, $\|\phi(y)-id(y)\|_0=\|r\psi(y/r)-r(y/r)\|\leq
%Cr\|\phi-id\|_0$.
%
%Also $\|D\phi(y).v-D(id)(y).v\|_0=\|D\psi(y/r).v-D(id)(y/r).v\|_0$
%
%Finally,
%$\|D\phi(x)-D\phi(x')/(x-x')^{\al}\|=\|D\psi(x/r)-D\psi(x'/r)/(r^{\al}(x/r-s'/r)^{\al}\|\leq
%C\|g-1\|_{\al}/r^{\al}$.
%
%And now, $\|g-1\|_{\al}=r^{\al}\|f-1\|_{\al}.$

\begin{remark}
A consequence of the Dacorogna-Moser's theorem is the fact that $\submundo$ is
path connected.
\end{remark}

So we obtain a weak version of the pasting lemma that allows us to change the diffeomorphism by its
derivative:

\begin{theorem}[Weak Pasting Lemma]\label{l.C2pasting}Let $M^n$ a compact Riemannian manifold without
boundary with dimension $n\geq 2$. If $f$ is a $C^2$-conservative
diffeomorphism and $x$ is a point in $M$, then for any $0<\alpha
<1,\eps >0$, there exists a $\eps$-$C^{1}$-perturbation $g$ (which
is a $C^{1+\alpha}$ diffeomorphism) of $f$ such that, for some small
neighborhoods $U\supset V$ of $x$, $g|_{U^c}=f$ and $g|_V=Df(x)$ (in
local charts).
\end{theorem}

\begin{proof} First consider a perturbation $h$ of $f$ such that
$h(y)=\rho(y)(f(x)+Df(x)(y-x))+(1-\rho(y))f(y)$ (in local charts),
where $\rho$ is a bump function such that $\rho|_{B(x,r/2)}=1, \
\rho|_{M-B(x,r)}=0, |\nabla\rho|\leq C/r\text{ and
}|\nabla^2\rho|\leq C/r^2$. Now, we note that $||h-f||_{C^1}\leq
C\cdot ||f||_{C^2}\cdot r$ and $||h-f||_{C^2}\leq
C\cdot||f||_{C^2}$, where $C$ is a constant. Then if we denote by
$\theta(y)=\det Dh(y)$ the density function of $h$, by the previous
calculation, we obtain that $\theta$ is $C\cdot ||f||_{C^2}\cdot
r$-$C^0$-close to $1$ (the density function of $f$) and $||\theta
-1||_{C^1}\leq C\cdot ||f||_{C^2}$. By the classical property of
convexity of the H\"older norms, we obtain that $||\theta
-1||_{\alpha}\leq C ||f||_{C^2}\cdot r^{1-\alpha}$. Hence,
increasing the constant $C$ if necessarily, we obtain the same bound
for the function $\widehat{\theta}(y)=(\det Dh(h^{-1}(y)))^{-1}$.
So, applying the lemma~\ref{l.DM4} for the domain
$\Omega=\overline{B(x,r)}-B(x,r/2)$ and the data function
$\widehat{\theta}$, we obtain that there exists a $C^{1+\alpha}$
diffeomorphism $\chi$ that is a solution to the equation of
lemma~\ref{l.DM4} which is $C\cdot ||f||_{C^2}\cdot
r^{1-\alpha}$-close to the identity. and is regular at the boundary
of $\Omega$ where $\chi=id$. Now observe that $g=\chi\circ h$ is
conservative and is close to $f$ provided that $r$ is small.
Smallness is obvious, and since $g(y)=f(x)+Df(x)(y-x)$ (in local
charts) in $B(x,r/2)$ we have that $\det Dg(y)=1$ there, and in
$M-B(x,r)$ we have $g(y)=f(y)$ so also $\det Dg(y)=1$, finally in
$\Omega$ we have that $\det Dg(x)=\det D\chi(h(x))Dh(x)=\det
D\chi(y)Dh(h^{-1}(y))=1$, so $g$ is conservative and this completes
the proof.
\end{proof}

\begin{remark}
We observe that the result also holds for $C^{1+\beta}$-diffeomorphisms, instead
of
$C^2$-diffeomorphisms. But in this case the perturbation is a
$C^{1+\alpha}$- diffeomorphism. Indeed, we can use the following interpolation
formula, which can be found in~\cite{Z}:
$$\|f\|_{(1-\gamma)a+\gamma b}\leq C\|f\|_a^{1-\gamma}\cdot\|f\|_b^{\gamma}.$$

\end{remark}

\begin{remark} In section~\ref{s.dom} of this paper, we apply this weak pasting
lemma to prove that every robust transitive conservative diffeomorphism admits a
dominated splitting. We hope that this lemma can be applied to obtain some other
interesting consequences. In fact, recently using the pasting lemma, Bochi-Fayad-Pujals were able to show
that $C^1$-stably $C^2$ ergodic systems are (generically) non-uniformly hyperbolic
(see~\cite{BFP}).
\end{remark}

\subsection{A Weak Pasting Lemma for Symplectic Diffeomorphisms}

The above perturbation lemma can be done for symplectic diffeomorphisms without
loosing the structure outside the neighborhood of a periodic orbit, as
follows. We also have the same result as in the conservative vector fields case,
if we require more differentiability of the diffeomorphism .

\begin{lemma}
If $f$ is a $C^k$-symplectic
diffeomorphism ($k\geq 1$), $x\in Per^n(f)$ is a periodic point of $f$ and $g$
is a local diffeomorphism ($C^k$-close to $f$) defined in a small
neighborhood $U$ of
the $f$-orbit of $x$,
then there exists a $C^k$-symplectic diffeomorphism $h$ ($C^k$-close to $f$)
and some neighborhood $U\subset V$ of $x$ satisfying $h|_U = g$ and
$g|_{V^c}=f$
\end{lemma}

\begin{proof} Consider the $f$-orbit $\SO(x)$ of $x$. Since all the perturbations are
local, by Darboux's theorem(see page 221 of~\cite{KH}), we can use
local coordinates near each point in $\SO(x)$, say
$\vfi_i:U_i\rightarrow V_i$, where $f^i(x)\in U_i$ and $0\in V_i =
B(0,\ga_i)\subset\real^n$. With respect to these coordinates, we
have the local maps $f:B(0,\de_i)\rightarrow V_{i+1}$, where
$V_n=V_1$, $B(0,\de_i)\subset V_i$ are small neighborhoods of $0$,
$f(0)=0$. In the symplectic case, we fix some bump function
$\la:\real\rightarrow [0,1]$ such that $\la(z)=1$ for $z\leq 1/2$
and $\la(z)=0$ for $z\geq 1$. Let $S_{1,i}$ be a \emph{generating
function} for $g:B(0,\de_i)\rightarrow V_{i+1}$ and $S_{0,i}$ a
generating function for $f:B(0,\de_i)\rightarrow V_{i+1}$
(see~\cite{X} for definitions and properties of generating
functions). Finally, set
$$S(x,y)=\la(\frac{2||(x,y)||}{\de_i})\cdot S_{1,i}(x,y) + \big[
1-\la(\frac{2||(x,y)||}{\de_i}) \big]\cdot S_{0,i}(x,y).$$ In
particular, if $h$ is the symplectic map associated to $S$, then
$h=g$ in $B(0,\de_i/4)$ and $g=f$ outside of $B(0,\de_i/2)$. To
summarize, this construction give us a symplectic
$\ep$-$C^1$-perturbation $g$ of $f$, if the numbers $\de_i$ are
small, such that $h$ coincides with $g$, near to each point
$f^i(x)$.
\end{proof}

\begin{remark} We can not obtain global versions of the pasting lemma for
symplectic diffeomorphism using this arguments, as proved for conservative ones,
since we used \emph{generating functions} which is a
local tool.
\end{remark}

In the next sections we study some other consequences of the pasting lemmas for
conservative systems.

%%%%%%%%%%%%%%%%%%%%%%%%%%%%%%%%%%%%%%%%%%%%%%%%%%%%%%%%%%%%%5
\section{Robustly Transitive Conservative Flows in 3-manifolds are Anosov}
\label{mpp}
Let $M^3$ be a compact three manifold and $\om$ a smooth volume in $M$. A
vector field $X\in \Mundo$ is $C^1$-\emph{robustly transitive} (in
$\Mundo$) if there is $\ep>0$ such that every $\ep$-$C^1$-perturbation $Y\in\Mundo$ of $X$ is transitive.

In this section we study conservative flows in 3-dimensional manifolds. Recently
the following result has been obtained by Morales, Pac\'{\i}fico and Pujals~\cite{MPP}: {\it Any isolated singular
$C^1$ robustly transitive set is a proper attractor.} This result is related
with a theorem by C. Doering~\cite{D}: {\it A $C^1$ robustly transitive flow is
Anosov}. In the same spirit of the diffeomorphism case, we can use the
pasting lemma to prove the conservative version of these results. We now
state the theorems:

\begin{maintheorem}
\label{t.doering}
Let $X\in \Mundo$ be a conservative $C^1$-robustly transitive vector field in a
3-dimensional compact manifold. Then $X$ is Anosov.
\end{maintheorem}

\begin{maintheorem}
\label{t.e} If $\Lambda$ is an isolated robustly transitive set of a
conservative vector field $X\in \Mundo$ in a 3-dimensional compact
manifold then it cannot have a singularity.
\end{maintheorem}

In particular, because the geometrical Lorenz sets are robustly transitive and
carrying singularities~\cite{Max}, as a corollary we have:

\begin{corollary}
There are no geometrical Lorenz sets (in the conservative setting)
for $C^1$ conservative vector fields in 3-dimensional compact
manifolds.
\end{corollary}

\subsection{Some Lemmas}

As in the diffeomorphism case we cannot have elliptic periodic orbits (or
singularities) in the presence of robust transitivity. In fact, generically we
have that any \emph{singularity} is hyperbolic by an argument of
Robinson~\cite{Rob}, but
in the case where we have robust transitivity, we can give another proof using
the pasting lemma.
More precisely we have the following lemmas:

\begin{lemma}
\label{l.noeliptic}
If $X\in \Mundo$ is $C^1$-robustly transitive then there are no elliptic
singularities (i.e. the spectrum of $DX(p)$ does not intersect $S^1$).
\end{lemma}

\begin{proof}
Suppose that there exist an elliptic singularity $\sigma$. Then we get $K$ a
small compact neighborhood of $\sigma$ such that $X_K$ is $C^1$-close to the linear
vector field $DX(\sigma)$. Now, we use the theorem~\ref{t.c1glue} to obtain a
vector field $Y$ $C^1$-close to $X$ which is $DX(\sigma)$ (in local charts)
inside an (possibly smaller) invariant
neighborhood. This contradicts transitivity.
\end{proof}

The novelty is that in the robustly transitive case, we also have the same
result for periodic orbits:

\begin{lemma}
\label{l.noelipticvector}
If $X\in \Mundo$ is $C^1$-robustly transitive then every periodic orbit is
hyperbolic.
\end{lemma}

\begin{proof}
Suppose that there exist an elliptic periodic orbit $p$. Then we get $K$ a
small compact neighborhood of $O(p)$ (e.g., a tubular neighborhood) such that
the linear flow induced by the periodic orbit restricted to $K$ is $C^1$-close
to $X$.
Now, we use the theorem~\ref{t.c1glue} to obtain a
vector field $Y$ $C^1$-close to $X$ with a (possibly smaller) invariant
neighborhood. This contradicts transitivity.
\end{proof}

\begin{remark}
\label{r.lemapraconjunto}
The above lemmas hold for robustly transitive sets strictly contained in
the whole
manifold with only minor modifications on the statement.
\end{remark}

\subsection{Proof of Theorem~\ref{t.doering}}

We already know that the singularities are hyperbolic.
Also by lemma~\ref{l.noelipticvector} every periodic orbit is hyperbolic. As
usual we will denote by $Crit(X)$ the set of periodic orbits and singularities.

\begin{remark} We observe that the hyperbolicity of hyperbolic periodic orbits
needs
the robust transitivity and the pasting lemma since Kupka-Smale's theorem
is false on dimension three
for conservative vector fields.
\end{remark}

Now we note that the result below by Doering~\cite{D} also holds in the conservative
case.
\begin{theorem}
\label{t.doeringcoisado} Let $\SS(M^3)=\{X\in\Mundo; \text{ every
}\sigma\in Crit(X) \text{ is hyper}$ $\textrm{bolic}\}$. If
$X\in\SS(M^3)$ then $X$ is Anosov.
\end{theorem}

Now if $X$ is robustly transitive then we already know that any critical element
is hyperbolic so we have the result.

We only sketch Doering's proof, since is straight forward that it works in the
conservative case:
\begin{lemma}
\label{l.doe} If $X\in\SS(M^3)$ then $M^3-Sing(X)$ has a dominated
splitting. In particular there are no singularities:
$Sing(X)=\emptyset$.
\end{lemma}

\begin{proof}

We start with a result on~\cite[Proposition 3.5]{D}:

{\bf Claim:} \textit{There exist $c>0$ and $0<\lambda<1$ and $\SU$ neighborhood of $X$ in
$\SS(M)$ such that there exist $(C,\lambda)$-dominated splitting for any $q\in
Per(Y)$ and $Y\in \SU$} (for the definition of $(C,\lambda)$ dominated splitting
for vector fields see~\cite{D}).

Let $x$ be a regular point for any $Y\in\SU$, so there are $Y_n$
such that $x\in Per(Y_n)$ (Pugh's Closing Lemma). Now we get the
$P^{Y_n}$-hyperbolic splitting over $O_{Y_n}(x)$ and by compactness
of the Grassmannian we have a splitting over $O_X(x)$ (by
saturation) then by the Claim we have that this splitting is a
$(C,\lambda)$-dominated splitting. Now by the abstract lemma of
Doering~\cite[Lemma 3.6]{D}:
\begin{lemma}
If $\Lambda$ is an invariant set of regular points of $X$ with a dominated
splitting over every orbit and such that every singularity of $X$ in the closure of
$\Lambda$ is a hyperbolic saddle then $\Lambda$ has a dominated splitting.
\end{lemma}
So using the lemma, we obtain a dominated splitting on $M-Sing(x)$.

By hyperbolicity $x_0$ is an interior point of $(M-Sing(X))\cup x_0$, but this
is a contradiction with the domination by the abstract~\cite[Theorem 2.1]{D}:

\begin{theorem}
Let $\Lambda$ an invariant set of regular points of $X$ and suppose that $x_0$
is a hyperbolic saddle-type singularity of $X$. If $\Lambda$ has a
dominated splitting then $x_o$ is not an interior point of $\Lambda\cup\{x_0\}$.
\end{theorem}

\end{proof}

The proof finish with the following result that is also available in
the conservative setting (see~\cite{Li1},~\cite{T1}, also
see~\cite{BDV} and references therein).

\begin{lemma}[Liao's theorem]
\label{l.liao} If $X\in\SS(M)$ without singularities then $X$ is
Anosov.
\end{lemma}

\subsection{Proof of the Theorem~\ref{t.e}}

We follow the steps of~\cite{MPP}. We denote by $\SU$ a neighborhood such that
every $Y\in\SU$ is transitive.

As in the proof of the previous theorem, we can suppose that any singularity or
periodic orbit is
hyperbolic. The first step is to prove that in fact it is a Lorenz-like
singularity. Recall that a singularity $\sigma$ is \emph{Lorenz-like} if its eigenvalues
$\lambda_2(\sigma)\leq\lambda_3(\sigma)\leq\lambda_1(\sigma)$ satisfy either $\lambda_3(\sigma)<0$ $\Rightarrow$
$-\lambda_3(\sigma)<\lambda_1(\sigma)$ or $\lambda_3(\sigma)>0$ $\Rightarrow$
$-\lambda_3(\sigma)>\lambda_2(\sigma)$.

\begin{lemma}
\label{l.lorenzlike} Let $M$ a three-dimensional compact manifold
boundaryless.   If $X\in \Mundo$ is $C^1$-robustly transitive
conservative vector field then any singularity $\sigma$ is
Lorenz-like.
\end{lemma}

\begin{proof}
First of all, the eigenvalues of $\sigma$ are real. Indeed, if $\omega=a+ib$ is an
eigenvalue of $\sigma$, then the others are $a-ib$ and $-2a$. Since the
singularity is hyperbolic, we can suppose also that $X\in \Submundo$ (using
theorem~\ref{l.dense}). By the Connecting
Lemma~\cite{WX} and~\cite{Ha}, we can assume that there exist a loop $\Gamma$ associated to
$\sigma$ which is a Shilnikov Bifurcation. Then by~\cite[page 338]{BS1} there
is a vector field $C^1$-close to $X$ with an elliptic singularity which gives a
contradiction with Lemma~\ref{l.noeliptic}(We observe
that the regularity for the bifurcation in~\cite{BS1} is more than or
equal to seven).

Let $\lambda_2(\sigma)\leq\lambda_3(\sigma)\leq\lambda_1(\sigma)$ the
eigenvalues. Now, $\lambda_2(\sigma)<0$ and $\lambda_1(\sigma)>0$, because
$\sigma$ is
hyperbolic and there are no sources or sinks ($\text{div }X=0$). Also, using that
$\sum \lambda_i =0$ we have that:
\begin{itemize}
\item $\lambda_3(\sigma)<0$ $\Rightarrow$
$-\lambda_3(\sigma)<\lambda_1(\sigma)$,
\item $\lambda_3(\sigma)>0$ $\Rightarrow$
$-\lambda_3(\sigma)>\lambda_2(\sigma)$.
\end{itemize}
\end{proof}

Now we give a sufficient condition which guarantees that $\Lambda$ is the whole
manifold.

\begin{theorem}
\label{t.nosing}
If $\Lambda$ is a transitive isolated set of a conservative vector field $X$,
such that:
\begin{enumerate}
\item $\Lambda$ contains robustly the unstable manifold of a critical element
$x_0\in Crit_X(\Lambda)$;
\item Every $x\in Crit_X(\Lambda)$ is hyperbolic;
\end{enumerate}
Then $\Lambda=M$.
\end{theorem}

\begin{proof}
We will use the following abstract lemma, which can be found in~\cite{MPP}:

\begin{lemma}[2.7 of~\cite{MPP}]
\label{l.27}
If $\Lambda$ is an isolated set of $X$ with $U$ an isolating block and a
neighborhood $W$ of $\Lambda$ such that $X_t(W)\subset U$ for every $t\geq 0$,
then $\Lambda$ is an attracting set of $X$.
\end{lemma}

We will prove that the hypothesis of the previous lemma are satisfied. This will
imply that $\Lambda=M$ since there are no attractors for conservative systems.
If there are no such $W$, then there exists $x_n\to x\in\Lambda$ and $t_n>0$
such that $X_{t_n}(x_n)\in M-U$ and $X_{t_n}(x_n)\to q\in \overline{M-U}$. Now,
let
$V\subset \ov{V}\subset int(U)\subset U$ a neighborhood of $\Lambda$. We have
$q\notin \ov{V}$.

There exists a neighborhood $\SU_0\subset \SU$ of $X$ such that
$\Lambda_Y(U)=\bigcap\limits_{t\geq0}Y_t(U)\subset V$ for every
$Y\in \SU_0$. And by hypothesis we have $W^u_Y(x_0(Y))\subset V$ for
a critical element $x_0\in Crit_Y(\Lambda_Y)$.

If $x\notin Crit_X(\Lambda)$, then let $z$ be such that $\Lambda=\omega_X(z)$
and $p\in
W^u_X(x_0)-O_X(x_0)$, so $p\in \Lambda$. Now we have $z_n\in O_X(z)$, $t_n'>0$
such that $z_n\to p$ and $X_{t_n'}(z_n)\to x$. Now, by the Connecting
Lemma~\cite{WX} and~\cite{Ha} there exist a $Z\in\SU_0$ such that $q\in V\supset W^u_Z(x_0(Z))$ and
$q\notin \ov{V}$, a contradiction.

If $x\in Crit_X(\Lambda)$ we can use the Hartman-Grobmann theorem to
find $x_n'$ in the positive orbit of $x_n$ ,$r\in
(W^u_X(x)-O_X(x))-Crit(X)$ and $t_n'>0$ such that $x_n'\to r$ and
$X_{t_n'}(x_n')\to q$. This is a reduction to the first case since
the following lemma says that $r\in \Lambda$  (the proof of the
lemma below is the same as in conservative case, since it only uses
the Connecting Lemma and the $\lambda$-lemma):

\begin{lemma}[2.8 of~\cite{MPP}]
\label{l.28}
If $\Lambda$ satisfies the hypothesis of theorem~\ref{t.nosing}, then
$W^u_X(x)\subset \Lambda$ for any $x\in Crit_X(\Lambda)$.
\end{lemma}
\end{proof}

Now we prove that there are no singularities. If $\Lambda=M$ then by
theorem~\ref{t.doering} $X$ will be Anosov, hence without singularities. So we
can assume that $\Lambda$ is a proper subset with a singularity and find a
contradiction. Now we follow the argument on~\cite{MPP}.

Taking $X$ or $-X$ there exists a singularity $\sigma$ such that
$\dim(W^u_X(\sigma))=1$(recall that dim$M$=3 and $\sigma$ is Lorenz-like). Let $U$ be an
isolating block such that $\Lambda_Y(U)$ is a connected transitive set for any
$Y\in\SU$. We will prove that $W^u_Y(\sigma(Y))\subset U$ for any $Y\in\SU$, so
by theorem~\ref{t.nosing} we have a contradiction.

Suppose that this does not happen for an $Y\in\SU$. By the dimensional hypothesis we
have two branches $w^+$ and $w^-$ of $W^u_X(\sigma)-\{\sigma\}$. Let $q^+\in
w^+$ and $q^-\in w^-$ (and let also $q^{\pm}(Y)$ be the continuation of
$q^{\pm}$ for $Y\in \SU$). Because the
negative orbit of $q^{\pm}(Y)$ converges to $\sigma(Y)$ and the unstable
manifold is not contained in $U$ there exists $t>0$ such that either
$Y_t(q^+(Y))$ or $Y_t(q^-(Y))$ is not in $U$. We assume that we are in the first case. We
know that there exists a neighborhood $\SU'\subset \SU$ of $Y$ such that for any
$Z\in \SU'$ we have $Z_t(q^+(Z))\notin U$. Let us take $z\in \Lambda_Y(U)$ with dense
orbit on $\Lambda_Y(U)$. This implies that $q^-(Y)\in \omega_Y(z)$. So we can
find a sequence $z_n\to q^-(Y)$ in $O_Y(z)$ and $t_n>0$ such that
$Y_{t_n}(z_n)\to q$ for some $q\in W^u_Y(\sigma(Y))-\{\sigma(Y)\}$. We set
$p=q^-(Y)$.

By the Connecting Lemma, there exists $Z\in\SU'$ such that $q^-(Z)\in
W^u_Z(\sigma(Z))$ and using the same arguments of theorem~\ref{t.nosing} we can find
$Z'\in\SU'$ and $t'>0$ such that $Z'_{t'}(q^-(Z'))\notin U$. This shows that
$\sigma(Z')$ is isolated in $\Lambda_{Z'}(U)$. By connectedness we get that
$\Lambda_{Z'}(U)$ is trivial, a contradiction.

%%%%%%%%%%%%%%%%%%%%%%%%%%%%%%%%%%%%%%%%%%%%%%%%%%%%%%%%%%%%%%%%%%%%%%
\vspace{.2cm}
%Using the same proof of theorem~\ref{t.morales} we can obtain the same result
%for conservative vector fields but in dimension greater than 3. So we have that
%generically if $X$ has a finite number of homoclinic classes then it is
%transitive. Indeed, by Pugh's lemma we obtain $M=Sing(X)\cup
%\bigcup\limits_{i=1}^n H(p_i)$ so there are no singularities and only one homoclinic
%class, hence transitive. This result cannot be directly strength to obtain transitivity of
%generic conservative flows, because the connecting lemma of Bonatti,
%Crovisier~\cite{BC} is not available in this case. Furthermore, we believe that far away from elliptic points it can be
%proved that there is a dominated splitting (the weak Herman's conjecture for
%flows). This will need a version for flows of the
%Bonatti-Diaz-Pujals~\cite{BDP} theorem.

\section{Robustly Transitive Conservative Diffeomorphisms}
\label{s.dom}

Now we study the existence of a dominated splitting for conservative
diffeomorphisms. This is a topic of many researches and several
results on the dissipative case are well known. To see that
transitive systems play an important role, we have the following
result of~\cite{BDP}: \emph{every $C^1$ robustly transitive
diffeomorphism has a dominated splitting}. This theorem has been
preceded by several results in particular cases :

\begin{itemize}
\item Ma\~ne~\cite{Ma}: \emph{Every $C^1$-robustly transitive diffeomorphism on
a compact surface is an Anosov system, i.e. it has a hyperbolic (hence dominated) splitting.}
\item Diaz, Pujals, Ures~\cite{DPU}: \emph{There is an open, dense set of 3-dimensional $C^1$-robustly transitive diffeomorphisms admitting dominated splitting}.
\end{itemize}

Let $f$ a $C^{1+\alpha}$-diffeomorphism, we say that $f$ is $C^1$-robustly
transitive if for any $g$ a $C^{1+\alpha}$-diffeomorphim that is $C^1$-close to
$f$ we have that $g$ is transitive. In this setting, we use the perturbation lemmas to prove
the conservative version of the theorem by Bonatti, Diaz and Pujals~\cite{BDP} mentioned above:

\begin{maintheorem}\label{BDP conserv.}Let $f$ be a
$C^{1+\alpha}$-diffeomorphism which is $C^1$-robustly transitive
conservative on a compact manifold $M^n$ with $n\geq
2$. Then $f$ admits a nontrivial dominated splitting defined on the whole $M$.
\end{maintheorem}

In dimension 2, we can extend Ma\~n\'e's theorem for $C^1$
conservative robustly transitive systems.

\begin{theorem}\label{Mane} A $C^1$-robustly transitive conservative system
$f\in\mundo$ is an Anosov diffeomorphism, where M is a compact surface.
\end{theorem}

Also we remark that an immediate corollary of
theorem~\ref{BDP conserv.} is :

\begin{corollary}\label{Ali} Let $f$ be a $C^1$-stably ergodic diffeomorphism in
$\submundo$. Then $f$ admits a dominated splitting.
\end{corollary}

This corollary answers positively a question posed by Tahzibi in~\cite{T}.

\subsection{Proofs of Theorem~\ref{BDP conserv.} and Theorem~\ref{Mane}}

First we give some definitions. We denote by $Per^n(f)$ the periodic points of $f$ with period $n$ and
$Per^{1\leq n}(f)$ the set of periodic points of $f$ with period at most $n$.
If $p$ is a hyperbolic saddle of $f$, the homoclinic class of $p$, $H(p,f)$, is the
closure of the transverse intersections of the invariant manifolds of $p$.

\begin{definition}
A compact invariant set $\Lambda$ admits a $(l,\lambda)$-dominated splitting for $f$ if
there exists a decomposition $T_{\Lambda}M=F\oplus G$ and a number $\lambda <1 $
such that for every $x\in \Lambda$:
$$\|Df^l|_{F(x)}\|\|Df^{-l}|_{G(x)}\|<\lambda.$$
\end{definition}

As remarked in~\cite{BDP}, it is sufficient to show that a conservative robustly
transitive diffeomorphism $f$ can not have periodic points $x\in Per^n(f)$ whose
derivative $D_x f^n$ is the identity. In fact, this is a consequence of :

\begin{proposition}[Proposition 7.7 in~\cite{BDP}]\label{dd de rob transit}
Given any $K>0$ and $\eps>0$ there is $l(\eps,K)\in\N$ such that for
every conservative diffeomorphism $f$ with derivatives $Df$ and $Df^{-1}$ whose
norm is
bounded by $K$, and for any saddle $p$ of $f$ having a non-trivial homoclinic class
$H(p,f)$, one has that:
\begin{itemize}
\item Either the homoclinic class $H(p,f)$ admits an $(l,1/2)$-dominated
splitting,
\item or for every neighborhood $U$ of $H(p,f)$ and $k\in \N$ there are a
conservative diffeomorphism $g$ $\eps-C^1$-close to $f$ and $k$ periodic points
$x_i$ of $g$ arbitrarily close to $p$, whose orbits are contained in $U$, such
that the derivatives $Dg^{n_i}(x_i)$ are the identity (here $n_i$ is the period
of $x_i$).
\end{itemize}
\end{proposition}

\begin{remark}
\label{r.smooth}
We stress that this perturbation $g$ preserves the differentiability class of the
original diffeomorphism, but $g$ is only $C^1$-close to $f$.
\end{remark}

Using the $C^{1+\alpha}$ pasting lemma~\ref{l.C2pasting} we have immediately the following:

\begin{lemma}
\label{Xia}If $f$ is a $C^{1+\alpha}$-conservative diffeomorphism and $x\in
Per^n(f)$ is a periodic point of $f$ such that $Df^n(x)$ is the identity matrix
then for any $\ep>0$, there is a $\ep$-$C^{1}$-perturbation $g$ of $f$
and some neighborhood $U$ of $x$ satisfying $g^n|_U = id$. And $g$ is a volume
preserving $C^{1+\alpha}$-diffeomorphism.
\end{lemma}

We now prove the theorem~\ref{BDP conserv.}:

\begin{proof}
We fix a constant $K>1$ which bounds the norms of the derivative of $f$ and $f^{-1}$ and
$l$ given by proposition~\ref{dd de rob transit}.

Now we use the following lemma whose proof we postpone.

\begin{lemma}
\label{l.estendeu}
There exists a small $C^1$-neighborhood $\SV$ of $f$ and $\delta>0$ such that
if $g\in \SV$ and $\Lambda\subset M$ is a $g$-invariant $\delta$-dense set with
a $(l,\lambda)$-dominated splitting for $g$ then $M$ has an
$(l,\lambda_1)$-dominated splitting for $g$ (where $\lambda_1$ is still less
than 1).
\end{lemma}

We fix the neighborhood of $\SV$ and $\delta$ given by the previous lemma.
Now we take a $\delta$-dense periodic saddle $p$ using the Closing
lemma in the conservative case~\cite{PR} (and again we do not loose
differentiability in this perturbation) for some $g$ in $\SV$ $C^1$-close to $f$.
So by Lemma~\ref{Xia} the second option in Proposition~\ref{dd de rob transit}
cannot be true, because choosing $\eps>0$ such that all of these perturbations
are in the neighborhood of $f$ where we have transitivity, we can find an invariant proper
open set for an iterate of $f$, and this contradicts transitivity. So the
homoclinic class $H(p,f)$ has a $l$-dominated splitting $F\oplus G$.
Then by Lemma~\ref{l.estendeu} we can extend the splitting for the whole
manifold and this defines a global $(l,\lambda_1)$-dominated splitting.
Hence we find a sequence $g_n$ which converges to $f$ with global $(l,\lambda_1)$-dominated
splitting, so $f$ also has a global dominated splitting. The proof is
complete.
\end{proof}

Now we will prove the lemma, but first we recall the notion of angles between
subspaces:

\begin{definition}
Let $V$ a finite dimensional vector space and $E$ a subspace of $V$. If
$F\subset V$ is a
subspace with the same dimension of $E$ and such that $F\cap E^{\perp}=\{0\}$ then
we define the angle between $F$ and $E$ as $\angle (F,E)=\|L\|$, where $L:E\to
E^{\perp}$ is a linear map such that $F=graph(L)$.

\end{definition}

\begin{proof}[Proof of Lemma~\ref{l.estendeu}]
First, we deal with the case that $g=f$. Observe that for every $\eps_0>0$
there exists $\gamma_0>0$ such that if $d(x,y)<\gamma_0$ then $x$ and $y$ are in a
trivializing local chart and
$\|Df(x)-Df(y)\|<\eps_0$. Also there exists $\beta>0$ such that for
any $(l,\lambda)$-dominated splitting $T_{\Lambda} M=E\oplus F$ we
have $\angle (E,F) >\beta$.

\begin{claim}
Given $\eps_1>0$ there exists $\gamma_1>0$ such that for any
$\Lambda$ with $T_{\Lambda}M=E\oplus F$ an $(l,\lambda)$-dominated
splitting if $d(x,y)<\gamma_1$ then $dist(F(x),F(y))<\eps_1$ (here
we are using the a distance in the Grassmannian).
\end{claim}
\begin{proof}
There exists an $m=m(l)$ such that if $v=v_1+v_2 \in E(x)\oplus F(x)$ then
$\|Df^m.v_1\|\geq C\lambda^m\|Df^m(x)|_{F(x)}\|$ (and a similar inequality holds
for $v_2$ and $E(x)$) and such that $C\lambda^m>>\sin \beta$.

Now we take $\gamma_1$ such that $d(f^j(x),f^j(y))<\gamma_0$ for every $j=0,...,m$.

Now suppose that $dist(F(y),F(x))\geq \eps_1$. Take $v\in F(y)$ with
$\|v||=1$ and write $v=v_1+v_2 \in E(x)\oplus F(x)$ (because $x$ and $y$ are in
the same local chart). So $\|v_2\|\leq C(\beta)$ a small constant, and
$\|v_1\|\geq C_1 \sin \alpha$ (where $\alpha$ depends of $\eps_1$, since we are
supposing that $F(x)$ is far from $F(y)$). So we have:
\begin{eqnarray*}
\|Df^m(y).v\| & \geq & \|Df(y)^m.v_1\|-\|Df(y)^m.v_2\| \\
 & \simeq & \|Df(x)^m.v_1\|-\|Df(x)^m.v_2\| \\
 & \geq & C\lambda^m\|Df(x)^m|_{F(x)}\|.\|v_1\|-\|Df(x)|_{F(x)}\|\|v_2\| \\
 & >> & 2\|Df^m(x)|_{F(x)}\|.
\end{eqnarray*}

So we have that $m(Df^m|_{F(y)})>>\|Df^m|_{F(x)}\|$ (here $m(L)$ is the co-norm
of the matrix $L$).

Now, since we are supposing that $dist(F(x),F(y))>\eps_1$, we can take $w\in
F(x)$ with $\|w\|=1$ and
write $w=w_1+w_2\in E(y)\oplus F(y)$. So we can
perform similar calculations to find that $m(Df^m|_{F(x)})>>\|Df^m|_{F(y)}\|$.
And this gives a contradiction.
\end{proof}

We know that there exists $K>0$ such that for any $\Lambda$ with an
$(l,\lambda)$-dominated splitting then $dist(Df(x)(C(x)),C(f(x)))>K$
for any $x\in \Lambda$ and $Df(x)C(x)\subset C(f(x))$. The next
claim says how we can obtain invariance of the cone fields in a
neighborhood of $\Lambda$.

\begin{claim}
Let $TM=E\oplus F$ a splitting, such that in an invariant compact
set $\Lambda$, this is a $(l,\lambda)$-dominated splitting. So we
have cone fields $C^*(x)$ ($*=s,u$), such that in $\Lambda$, these
are invariant cone fields associated to the dominated splitting.
Suppose that there exists $\gamma$ such that, $d(x,y)<\gamma$
implies $dist(C^*(x),C^*(y))<K/4$ (in the projective metric). Then
there exists $\delta>0$ (which depends only of $(l,\lambda)$ and the
constants of uniform continuity of the diffeomorphism and the
derivative) such that in a $\eta$-neighborhood of $\Lambda$, the
cone fields are also invariant.
\end{claim}

\begin{proof}
We deal with $C^E(x)$, the case $C^F(x)$ is similar. We know that
$$dist(Df(x)(C(x)),C(f(x)))>K$$for any
$x\in \Lambda$ and $Df(x)C(x)\subset C(f(x))$ (here we are using a projective
distance on cones). By uniform
continuity, there exists $\delta_1$ such that if $d(x,y)<\delta_1$ then
$$dist(C(f(x)),C(f(y)))$$is small, then (trivializing the tangent bundle) we have
that $$dist(Df(x)(C(x)),C(f(y)))>K/2 \textrm{ and }Df(x)C(x)\subset
C(f(y)).$$Now, by hypothesis, we have that if
$d(x,y)<\gamma$,  then we have that
$$dist(Df(x)(C(y)), C(f(y)))>K/4$$and $Df(x)C(y)\subset C(f(y))$. Finally, by uniform continuity of $Df(x)$, we
have that if $d(x,y)<\delta_2$ we have that $dist(Df(y)C(y),C(f(y)))>K/8$ and
$Df(y)C(y)\subset C(f(y))$. This proves the claim.
\end{proof}

We will say that a function $H:M\to Y$ is $(\eps,\eta)$-continuous
on $X\subset M$ if for any $a,b\in X$ such that $d_X(a,b)<\eta$ then
$d_Y(f(a),f(b))<\eps$. Now by the first claim, we have that for any
$\Lambda$ with an $(l,\lambda)$-dominated splitting, the bundles are
$(\eps,\eta)$-continuous in $\Lambda$. Take $\eps$ such that if the
bundles are $(4\eps,\eta/2)$-continuous then the cone fields are
$(K/16,\eta)$-continuous.

So, if we have that $\Lambda$ is $\eta/2$-dense, then take $(U_i)$ a
covering of $\Lambda$ by balls of radii $\eta$, hence is also a
covering of $M$. So each bundle restricted to $U_i\cap\Lambda$ have
variation bounded by $2\eps$. So by Tietze's theorem we can extend
the bundles in the whole $U_i$, obtaining bundles $\widehat{E}$ and
$\widehat{F}$ on $U_i$ with variation bounded by $2\eps$, and so we
can glue this bundles on $M$ to obtain $\widehat{E}$ and
$\widehat{F}$ bundles which are $(4\eps,\eta/2)$-continuous in $M$.

By the second claim, using $\gamma=\eta/2$, the cone fields
associated to $\widehat{E}$ and $\widehat{F}$ are invariant in a
neighborhood of size $\delta$. But the previous analysis holds for
any $\Lambda$ with an $(l,\lambda)$-dominated splitting. So if we
take $\Lambda$ also $\delta$-dense. We can extend the splitting
$T_{\Lambda}M=E\oplus F$ to the whole manifold, and obtain invariant
cone fields in the whole manifold. This gives another splitting
$TM=\widetilde{E}\oplus\widetilde{F}$ which is dominated.

It's clear now, that the proof for $f$, uses only the constant of uniform
continuity of $f$ and $Df$ and also the norm $\|f\|_{C^1}$. But these are
uniforms in a $C^1$-neighborhood of $f$, then the lemma follows.
\end{proof}

Now we prove the theorem~\ref{Mane}:

\begin{proof}
We will use the following theorem, which holds in the $C^1$-topology for
$C^1$-diffeomorphisms:
\begin{theorem}[Theorem 6 of~\cite{BDP}]
\label{t.6bdp} Let $f$ be a conservative diffeomorphism. Then
either:
\begin{enumerate}
\item for any $k\in \N$ there
exists a conservative diffeomorphism $g$ arbitrarily $C^1$-close to $f$ having
$k$ periodic orbits whose derivatives are the identity or;
\item $M$ is the union of
finitely many invariant compact sets having a dominated splitting.
\end{enumerate}
If $f$ is
transitive and the second possibility occurs then $f$ admits a dominated
splitting.
\end{theorem}

Now to avoid periodic points whose derivative is the identity we use the
pasting lemma in the symplectic case since the dimension of the manifold is
two, and this can be done for $C^1$ diffeomorphisms. So we obtain a global
dominated splitting. To finish we use:

\begin{proposition}[Proposition 0.5 of~\cite{BDP}]
Let $f$ be a conservative diffeomorphism with a dominated splitting. Then the
derivative contracts uniformly the volume of one of the bundles and expands
uniformly the volume of the other bundle.
\end{proposition}

So because in dimension two, both subbundles are one-dimensional we have
hyperbolicity.
\end{proof}

\section{Appendix: On a non-estimate theorem for the divergence equation (by
David Diica and Yakov Simpson-Weller)}
\label{s.dica}

We show that a non-estimate result (similar to the well-known Ornstein's
non-estimate theorem) holds for the divergence equation. In particular, the
method of construction of solutions of the Jacobian determinant equation by
Dacorogna-Moser does not work in the limit case of $C^1$ diffeomorphisms.

\subsection{Introduction and Historical Notes}

The main result of this appendix is :
\begin{maintheorem}\label{t.1} There exist no operator (even \emph{nonlinear})
$$\SL:B_r^0\cap\newmundo\rightarrow\newmundo$$ such that $\text{div}\SL g=g$ and
$||\SL(g_1)-\SL(g_2)||_{C^1}\leq ||g_1-g_2||_{C^0}$ (i.e., $\SL$ is Lipschitz
with respect to the $C^0$ norm $||.||_{C^0}$ and the $C^1$ norm $||.||_{C^1}$).
\end{maintheorem}

The motivation of the previous theorem is the classical problem of denseness of
smooth conservative diffeomorphisms of compact manifolds. This problem was successfully
solved for $C^{1+\al}$ diffeomorphisms, using the regularity estimates for the
Laplacian operator. In fact, the problem of denseness of smooth diffeomorphism
can be reduced to the problem of solving the Jacobian determinant equation
$\det\nabla u = f$, which can be linearized to obtain the equation $\text{div} v
=g$. If we look for special solutions $v=\text{grad } a$ of this equation, the
problem is reduced to solve the equation $\Delta a =g$. The idea is then to
``deform'' (Moser's ``flow idea'')
the solutions of the equation $\Delta a=g$ to solutions of the Jacobian
determinant equation. But, the standard elliptic regularity theory only
guarantees solutions for H\"older functions. Using this ideas, Zehnder was able
to positively answers the denseness problem for $C^{1+\al}$ systems. See, for
example, Dacorogna-Moser's work~\cite{DM} (and the reference to Zehnder's article
therein) for more details.

However, the method used by Dacorogna-Moser and Zenhder to solve the divergence
equation (and, in particular, the determinant equation) is \emph{linear}. So, in
order to solve the denseness problem for $C^1$ conservative diffeomorphisms, one
can look for \emph{nonlinear} bounded right inverses for the divergence
operator. Recently, Bourgain and Brezis~\cite{BB} have shown that there are
\emph{explicit nonlinear} bounded right inverses for the divergence operator,
but there are no bounded linear right inverses. The problem to apply their
result is that their theorem is stated for the Sobolev spaces $W^{1,d}$.  On the
other hand, since Burago-Kleiner and McMullen (see~\cite{BK}) showed that for some
$L^{\infty}$ (and even, continuous) functions, the determinant equation admits
no solution. The idea of McMullen is that if the equation admits solutions for
every $L^{\infty}$ function, then we have a contradiction with a well-known
``non-estimate'' result due to Ornstein~\cite{O}.

In this short note, we consider a lemma due to Burago-Kleiner~\cite{BK} in order
to show that there are no Lipschitz estimates for the nonlinear right inverses
constructed by Bourgain-Brezis. In particular, Dacorogna-Moser's idea can not
work for the limit case of $C^1$ diffeomorphisms, even when we consider the
``nonlinear'' modification of the idea, which consists in substituting the
linear right inverse of the divergence
operator in the H\"older case by Lipschitz nonlinear right inverses.

\subsection{Main Tools}

In this section we state some useful facts which will be used in the proof of the
non-estimate theorem. First of all, we recall the following result due to
Burago-Kleiner (see~\cite[page 275]{BK}):

\begin{lemma}[Burago-Kleiner]\label{l.1}For every $L>1$, $c>0$, there exists a smooth
function $\rho_{L,c}: I^d\rightarrow [1, 1+c]$ (where $I^d$ is the
$d$-dimensional cube, $I=[0,1]$) such that the equation
$\det(\nabla\psi)=\rho_{L, c}$ does not admit an $L$-biLipschitz solution
$\psi$.
\end{lemma}

\begin{remark} As an immediate corollary of this lemma, it follows that the
following non-estimate result holds for the Jacobian determinant equation:
\emph{There are no bounded right inverses for the Jacobian determinant
equation, i.e., there exists no operator $\SL$ such that $|| \SL f||_{C^1}
\leq K\cdot ||f||_{C^0}$ for every $f\in C^{\infty}$ $C^0$-close to the constant
function $1$ and $\det\nabla\SL (f) =f$}.
\end{remark}

Now we briefly recall the method (by Dacorogna-Moser~\cite{DM}) of construction of solutions of the Jacobian
determinant equation $$\det \nabla u = f  $$ using solutions of the
divergence equation $$\text{div} Y= f .$$

Let for any $n\times n$ matrix $\zeta$, $$Q(\zeta)=\det (I+\zeta)-1-\text{tr}
(\zeta) ,$$ where $I$ denotes the identity matrix. Since $Q$ is a sum of
monomials of degree $t$ with $2\leq t\leq n$, we can find a constant $K_2 >0$
such that if $w_1, w_2\in C^k(\R^n,\R^n)$ with $||w_1||_0 , ||w_2||_0 \leq 1$, then $$ (*)
|| Q(w_1)-Q(w_2)||_k\leq K_2 (||w_1||_0 + ||w_2||_0) ||w_1 -w_2 ||_k .$$ Here
$k$ is any real number, but the case of interest for this paper is $k$ integer.

With this notation, Dacorogna-Moser's idea to construct solutions of the
Jacobian determinant equation is: denote by $u(x)=x+v(x)$ and \emph{suppose}
that there is an operator $\SL : \mundo\rightarrow\submundo$ (\emph{not
necessarily linear}) such that $\text{div}\SL(a)=a$, for every $a\in\mundo$, and
$\SL$ is Lipschitz continuous. For example, there are bounded linear operators
$\SL$ as above when $k=r +\alpha$, $r\geq 0$ integer, $\alpha >0$. Define the
operator $N:\submundo\rightarrow\mundo$, $$N(v)= f-1-Q(\nabla v) .$$ Clearly
any fixed point $v$ of the map $\SL\circ N$ (i.e., $v=\SL N(v)$) satisfies
$\det\nabla u=f$ (recall that $u=\text{id}+v$). So to obtain solutions of the
Jacobian determinant equation, it suffices to show that, if $||f-1||_k$ is
small, then $\SL\circ N$ is a contraction. But, it is not difficult to see that
this follows from the inequality $(*)$ and $||f-1||_k$ is small
(see~\cite[page 12]{DM} for more details).

\begin{remark}If $||f-1||_k$ is not small, Dacorogna-Moser uses a trick
which consists of two steps: First, considering the quotient of $f$ by a
suitable smooth function $g$, the quotient has small norm. Using the ``flow
idea'', the case of the determinant equation for a smooth density $g$ is
solved (see~\cite[page 13]{DM}).
\end{remark}

After these considerations, we are in position to prove the main theorem of this
appendix.

\subsection{Proof of theorem~\ref{t.1}}

Our proof of theorem~\ref{t.1} is by contradiction.

\begin{proof}[Proof of theorem~\ref{t.1}] Suppose that there exists some
Lipschitz continuous (\emph{not necessarily linear}) operator $\SL$ such that
$\text{div}\SL g=g$ for every $f\in B_r^0\cap\newmundo$ ($B_r^k$ is the $r$-ball in
the space of $C^k$ functions with the $C^k$-norm). Fix some positive
number $L>0$ such that $L>K_1\cdot r$, where $K_1$ is the Lipschitz constant of
$\SL$. Then, if we consider the smooth function $\rho_{L,c}$ given by
Lemma~\ref{l.1}, where $c$ is sufficiently small, then, by the facts remarked in
the previous section, $\SL\circ N$ is a contraction ($N(v)= \rho_{L,c}-1-
Q(\nabla v)$) and, in particular, there
exists some fixed point $v$ of $\SL\circ
N : (B^1_{\epsilon}\cap\newmundo)\rightarrow (B^1_{\epsilon}\cap\newmundo)$
($\eps$ is chosen such that
$N : B^1_{\epsilon}\cap\newmundo\rightarrow B_r^0\cap\newmundo$; $\eps>0$ exists by
the inequality $(*)$). Moreover, $u=\text{id}+v$ solves
the equation $\det\nabla u=\rho_{L,c}$, and the $C^1$-norm of $u$ is at most
$K_1\cdot r$. In particular, since $K_1\cdot r<L$, $u$ is $L$-biLipschitz (for
$r$ small). But this contradicts the conclusion of lemma~\ref{l.1}.
\end{proof}

\begin{ack}
The authors are thankful to Sylvain Crovisier, Christian Bonatti, Vitor Ara\'ujo
and Maria Jos\'e Pac\'{\i}fico for useful
conversations and remarks. We also are indebt to Jairo Bochi who found a gap in
a previous version of the paper and for his comments. We want to thank also Leonardo
Macarini who pointed out to our attention Dacorogna-Moser's theorem. The
authors also want to stress their gratitude to Marcelo Viana for his
encouragement and guidance. Finally, we want to thank IMPA and its staff.

\end{ack}

%%%%%%%%%%%%%%%%%%%%%%%%%%%%%%%%%%%%%%%%%%%%%%%%%%%%%%%%%%%%%%

\bibliographystyle{alpha}
\bibliography{bib}

\begin{thebibliography}{MPP}

\bibitem[Au]{Au}
T.~Aubin,
\newblock Nonlinear Analysis on Manifolds. Monge-Amp\`ere Equations,
\newblock Springer-Verlag, 1982.

\bibitem[BB]{BB}
J. Bourgain and H. Brezis,
\newblock On the equation $\text{div} Y=f$ and applications to control of
phases,
\newblock \emph{J. of the A. M. S.}, vol.16 (2), 393--426, 2002.

\bibitem[BK]{BK}
D. Burago and B. Kleiner,
\newblock Separated nets in Euclidean spaces and Jacobians of biLipschitz maps,
\newblock \emph{Geometric and Funct. Anal.}, vol.8, 273--282, 1998.

\bibitem[BFP]{BFP}
J.~Bochi, B.~Fayad and E.~Pujals.
\newblock A remark on conservative diffeomorphisms.
\newblock \ \ In preparation.

%\bibitem[BC]{BC}
%C. Bonatti and S. Crovisier,
%\newblock R\'ecurrence et G\'en\'ericit\'e,
%\newblock Preprint 2003.

\bibitem[BDP]{BDP}
C.~Bonatti, L.~Diaz, E.~Pujals,
\newblock A $C^1$-generic dichotomy for diffeomorphisms: Weak forms of
hyperbolicity or infinitely many sinks or sources,
\newblock Preprint 2002, to appear in \emph{Annals of Math}.

\bibitem[BDV]{BDV}
C.~Bonatti, L.~Diaz, M.~Viana,
\newblock Dynamics Beyond Uniform Hyperbolicity,
\newblock Springer Verlag, 2005.

\bibitem[BS]{BS1}
V. Biragov and L. Shilnikov,
\newblock On the bifurcation of a saddle-focus separatrix loop in a
three-Dimensional conservative dynamical system,
\newblock \emph{Selecta Mathematica Sovietica}, v.11, no.4, 333--340, 1992.

%\bibitem[BV]{Bochi-Viana}
%J.~Bochi and M.~Viana,
%\newblock The Lyapounov exponents of generic volume preserving and symplectic
%systems,
%\newblock Preprint 2002, to appear in \emph{Annals of Math.}

%\bibitem[CMP]{CMP}
%C.~Carballo, C.~Morales and M.~Pac\'{i}fico,
%\newblock Homoclinic classes for generic $C^1$ vector fields,
%\newblock Preprint 2001, to appear in \emph{Erg. Theory and Dyn. Systems}.

\bibitem[DM]{DM}
B.~Dacorogna, J.~Moser,
\newblock On a partial differential equation involving the jacobian determinant,
\newblock \emph{Ann. Inst. Poincar\'e}, V.7. 1--26. 1990.

\bibitem[DPU]{DPU}
L.~Diaz, E.~Pujals, R~Ures,
\newblock Partial hyperbolic and robust transitivity.
\newblock \emph{Acta. Math.}, v. 183, 1--42, 1999.

\bibitem[D]{D}
C.~Doering.
\newblock Persistently transitive vector fields on three-dimensional manifolds,
\newblock \emph{Proc. on Dynamical Systems and Bifurcation Theory}, Pitman Res.
Notes Math. Ser., 160, 59--89, 1987.

\bibitem[GT]{GT}
D.~Gilbarg and N.~Trudinger
\newblock Elliptic Partial Differential Equations of Second Order,
\newblock Springer-Verlag, 1998.

%\bibitem[Her]{Her}
%M.~Herman,
%\newblock Some open problems in dynamical systems.
%\newblock \emph{Doc. Mathematica. Extra Volume ICM}, 1998-II, 797--808.

\bibitem[Ha]{Ha}
S.~Hayashi,
\newblock Connecting invariant manifolds and the solution of the $C^1$ stability
and $\Omega$-stability conjectures for flows.
\newblock \emph{Annals of Math.}, 145, 81-137, 1997.

\bibitem[H]{H}
L.~H\"ormander,
\newblock Linear Partial Differential Operators,
\newblock Springer-Verlag, 1976.

\bibitem[KH]{KH}
A.~Katok and B.~Hasselblat,
\newblock Introduction to the Modern Theory of Dynamical Systems.
\newblock Cambridge University Press, 1995.

\bibitem[Li1]{Li1}
S. T. Liao,
\newblock Obstruction sets II,
\newblock \emph{Beijin Daxue Xuebao},no 2, 1--36, 1981.

\bibitem[LU]{LU}
O.~Ladyzenskaya and N.~Uralsteva,
\newblock Linear and Quasilinear Elliptic Partial Differential Equations,
\newblock Academic Press, 1968.

\bibitem[Ma]{Ma}
R.~Ma\~n\'e,
\newblock An ergodic closing lemma.
\newblock \emph{Ann. Math.}, v. 116, 503--540. 1982.

\bibitem[MPP]{MPP}
C.A. Morales, M.J. Pac\'{\i}fico, E. Pujals.
\newblock Robust transitive singular sets for 3-flows are partially hyperbolic
attractors or repellers.
\newblock To appear in the Annals of Mathematics.

\bibitem[Mo]{M}
J.~Moser,
\newblock On the volume element of a manifold.
\newblock \emph{Trans. AMS}, v.120, 286--294, 1965.

%\bibitem[N]{N}
%S.~Newhouse,
%\newblock Quasi-elliptic periodic points in conservative dynamical systems.
%\newblock \emph{Amer. J. Math.}, V. 99, no. 5, 1061--1087, 1977.

\bibitem[O]{O}
D. Ornstein,
\newblock A non-inequality for differential operators in the $L^1$ norm,
\newblock \emph{Arch. Rat. Mech. and Anal.}, vol.11, 40--49, 1962.

\bibitem[PR]{PR}
C.~Pugh and C.~Robinson,
\newblock The $C^1$ Closing Lemma, including hamiltonians,
\newblock \emph{Ergodic Theory and Dynamical Systems}, vol. 3, no.2, 261--313,
1983.

\bibitem[R]{Rob}
R.C.~Robinson,
\newblock Generic Properties of conservative systems, I and II,
\newblock \emph{Amer. J. Math}, 92, 562--603 and 897--906, 1970.

\bibitem[RY]{RY}
T.~Riviere and D.~Ye,
\newblock Resolutions of the prescribed volume form equation,
\newblock \emph{Nonlinear Diff. Eq. Appl.}, 3, 323--369, 1996.

\bibitem[T]{T}
A.~Tahzibi,
\newblock Stably ergodic systems which are not partially hyperbolic,
\newblock Preprint 2002.

\bibitem[T1]{T1}
H.~Toyoshiba,
\newblock Nonsingular vector fields in $\SG^1(M^3)$ satisfy Axiom A and no cycle: a new proof of Liao's theorem,
\newblock \emph{Hokaido Math. J.}, 29, no. 1, 45--58. 2000.

\bibitem[V]{Max}
M.~Viana,
\newblock What's new on Lorenz Strange Attractors
\newblock \emph{Math. Inteligencer}, vol. 22, no. 3, 6--19, 2000.

\bibitem[WX]{WX}
L.~Wen and Z.~Xia,
\newblock $C^1$ connecting lemmas,
\newblock \emph{Trans. Amer. Math. Soc.}, 352, no. 11, 5213--5230, 2000.

\bibitem[X]{X}
Z.~Xia,
\newblock Homoclinic points in symplectic and volume preserving diffeomorphisms,
\newblock \emph{Comm. in Math. Phys.}, Vol.177 , 435--449, 1996.

\bibitem[Z]{Z}
E.~Zehnder,
\newblock Note on smoothing symplectic and volume preserving diffeomorphisms,
\newblock \emph{Lect. Notes in Math.}, 597, 828--854, 1977.

\bibitem[Zu]{Zu}
C.~Zuppa,
\newblock Regularisation $C^{\infty}$ des champs vectoriels qui pr\'eservent
l'el\'ement de volume.
\newblock \emph{Bol. Soc. Brasileira Matem.}, 10 (2), 51--56, 1979.

\end{thebibliography}

\vfill

%\vspace{1cm}
{\footnotesize
\noindent

\vfill

Alexander Arbieto ({\tt alexande{\@@}impa.br}) \\

Carlos Matheus ({\tt cmateus{\@@}impa.br}) \\

\smallskip

\noindent IMPA, Estrada D. Castorina 110, Jardim Bot\^anico, 22460-320 Rio de Janeiro, Brazil
}

\end{document}